\documentclass[reqno]{amsproc}
\usepackage[latin1]{inputenc}
\usepackage{amssymb}
\usepackage{amsmath}
\usepackage{relsize,exscale}
\usepackage{latexsym}%
\usepackage[frame,ps,matrix,arrow,curve,rotate,all,2cell,tips]{xy}
\usepackage{cite}
\usepackage{pgf,tikz}
\usetikzlibrary{matrix}
\usetikzlibrary{arrows}
\usepackage{amssymb}
\usepackage{amsfonts}
\usepackage{amscd}
\usepackage{xcolor}
\usepackage{hyperref}
\usepackage{amstext,amsmath,amssymb,amsfonts}
\newtheorem{theorem}{Theorem}

\newtheorem{corollary}[theorem]{Corollary}

\newtheorem{definition}[theorem]{Definition}

\newtheorem{lemma}[theorem]{Lemma}

\newtheorem{proposition}[theorem]{Proposition}
\newtheorem{remark}[theorem]{Remark}

\newcommand{\beq}{\begin{eqnarray}}
\newcommand{\eeq}{\end{eqnarray}}
\newcommand{\beqs}{\begin{eqnarray*}}
\newcommand{\eeqs}{\end{eqnarray*}}
\newcommand{\bpro}{\begin{pro}}
\newcommand{\epro}{\end{pro}}
\newcommand{\blem}{\begin{lem}}
\newcommand{\elem}{\end{lem}}
\newcommand{\bdfn}{\begin{dfn}}
\newcommand{\edfn}{\end{dfn}}
\newcommand{\bcor}{\begin{cor}}
\newcommand{\ecor}{\end{cor}}
\newcommand{\bthm}{\begin{thm}}
\newcommand{\ethm}{\end{thm}}
\newcommand{\bex}{\begin{ex}}
\newcommand{\eex}{\end{ex}}
\newcommand{\brmk}{\begin{rmk}}
\newcommand{\ermk}{\end{rmk}}
\newcommand{\bpr}{\begin{pr}}
\newcommand{\epr}{\end{pr}}
\newcommand{\benum}{\begin{enumerate}} 
\newcommand{\eenum}{\end{enumerate}}
\newcommand{\bitem}{\begin{itemize}}
\newcommand{\eitem}{\end{itemize}}

\chardef\bslash=`\\
\numberwithin{equation}{section}
\numberwithin{table}{section}
\numberwithin{theorem}{section}

\title[Fundamental Theorems of Fourier-Stieltjes Transform on Locally Compact Group]{Fundamental Theorems of Fourier-Stieltjes Transform Defined by Induced Representation on Locally Compact Group}
\author{Y. I. Akakpo $^\ast$}
\address[$\ast, \dagger, \ddagger$]{University of Abomey-Calavi,
International Chair in Mathematical Physics and Applications
(ICMPA-UNESCO Chair), 072 BP 50, Cotonou, Rep. of Benin}
\address[$ \ddagger$]{University of California, Los Angeles,
Department of Mathematics,
520 Portola Plaza, Los Angeles, CA 90095 }
\email{$^\ast$ akjeremiji@yahoo.fr}
\email{$^\dagger$ kofassiam@yahoo.fr}
\email{$^\dagger$ koffi@math.ucla.edu}
\author{ M. N.  Hounkonnou$^\ddagger$}
\author{K. Enakoutsa$^\dagger$}
\author{V. S. K. Assiamoua$^\dagger$}
\email{$^\ddagger$ norbert.hounkonnou@cipma.uac.bj, with copy to hounkonnou@yahoo.fr}
\begin{document}
\maketitle
\begin{abstract}
This work addresses an extension  of Fourier-Stieltjes transform of a vector measure  defined on compact groups to locally compact groups, by using a group representation induced by a representation of one of its compact subgroups. 
\end{abstract}



\noindent{
\today}
\tableofcontents

\section*{Introduction} 
The Fourier transform is one of the most important and widely used computational tasks. It is nowadays an essential tool in many technological innovations.

A foundational tool commonly used to analyze the spectral representation of signals. Its applications include audio/video processing, radar and GPS systems, wireless communications, medical imaging and spectroscopy, the processing of seismic data, and many other tasks. 

Among the Fourier transforms, we have the Fourier transform of integrable functions and that of measures.

During the last decades, countless studies on vector measures have been made \cite{Di}, \cite{Bar} (and references therein), and also \cite{assiamoua1} for some applications on compact groups. 


 The project tackled in this paper is a generalization of the work of Assiamoua \cite{assiamoua1, assiamoua2} on the Fourier-Stieltjes transformation of a vector measure and that of an absolutely continuous integrable function with respect to a Haar measure. This is a continuation of recent work by Akakpo and al \cite{yao} in which they set up the foundations of the Fourier-Stieltjes transform of a vector measure and that of an integrable function of a locally compact group $G$.

The paper is organized as follows. In Section 1,  we  give some preliminaries useful in the sequel. In Section 2, we construct some spaces called $\mathfrak{S}_{p}$ spaces, very useful in our theory, and give some of their main properties, show that the Fourier transform define in \cite{yao} verify   the main fundamental theorems in analysis. 

\section{Preliminaries and notations}
In this section, for the clarity of the development, we briefly recall main useful  known definitions and results, and set our notations.

We consider  a locally compact space $G,$ the  Banach spaces  $A$ and $F$  over the field $\mathbb{K}, \,  ( \mathbb{K} = \mathbb{R}$  or $\mathbb{C}),$
and denote by $\mathcal{K}(G, A)$ the vector space of all continuous functions $f : G \longrightarrow A$  with a compact support, and by $\mathfrak{C}(G,A)$ the space of continuous functions $f : G \longrightarrow A$.

As a matter of notation simplification, we write $\mathcal{K}(G)$ instead of $\mathcal{K}(G, \mathbb{R})$ or $\mathcal{K}(G, \mathbb{C}).$
 For each subset $K$ of $G$, denote by 
 $\mathcal{K}_{K}(G, A)$ the space  of functions with  support contained in $K$. $\mathcal{K}_{K}(G, A)$ is obviously a subspace of $\mathcal{K}(G, A)$.
\begin{definition} \mbox{ }\\ 
  For every function $f \in \mathcal{K}(G, A)$, we define  as usual
 \[\Vert f \Vert := \sup _{t\in G} \Vert f(t) \Vert_{A}. \] 
 The mapping $ f \mapsto \Vert f \Vert$ is a norm on each space  $\mathcal{K}_{K}(G, A),$ and defines  the topology of uniform convergence on $G$ over $\mathcal{K}(G)$. 
\end{definition}
\begin{definition} \mbox{ }\\ 
 On $\mathcal{K}(G, A),$  the topology of the  compact
convergence is the locally convex topology defined by the family of seminorms  \[\left\Vert f \right\Vert_{K} = \sup_{t \in K} \Vert f(t) \Vert_{A}, \]
where $K$ runs over the set of  compact subsets of $G$.
\end{definition}
\begin{proposition} \mbox{ }\\
 The space $\mathcal{K}(G, A)$ is dense in the space $\mathfrak{C}(G,A)$ for the
topology of the compact convergence.
\end{proposition}
\subsection{Vector measure}
\begin{definition} \mbox{ }\\ 
 A vector measure on $G$ with respect to two
spaces $A$ and $F$, or an $(A, F)$-measure on $G$, is any linear mapping 
$m : \mathcal{K}(G, A)\longrightarrow F $ having the property that,  for each compact set $K \subset G,$ the restriction $m$ to the subspace $\mathcal{K}_{K}(G, A)$ is continuous for the topology of uniform convergence, i.e for each compact set,  there exists a number $a_{K} > 0$ such that 
\[ \left\Vert m (f) \right\Vert \leq a_{K} \sup \left\lbrace \Vert f(t) \Vert_{A} , \quad t \in K \right\rbrace. \]

The value  $m(f)$ of $m$  for a function $f \in \mathcal{K}(G, A)$ is called the
integral of $f$ with respect to $m,$ also denoted by $\int_{G}f dm$ or $\int_{G}f(t) dm(t)$. 
A vector measure is said to be dominated if there exists a positive measure $\mu$ such that 
\[\left\Vert \int_{G}f(t)dm(t) \right\Vert \leq \int_{G} \vert f(t) \vert d\mu(t), \quad    \quad f \in \mathcal{K}(G). \]
If $m$ is dominated, then there exists a smallest positive measure $\vert m \vert,$ called the modulus or the variation of $m$, that dominates it.  A positive measure is said to be bounded if it is continous in the uniform norm topology of $\mathcal{K}(G).$
A vector measure is said to be bounded if it is dominated by a bounded positive measure. If $m$ is bounded,  then $\vert m \vert$ is also bounded. 

Denoting by $M_{1}(G,A)$ the Banach algebra of bounded vector measures on $G,$
 the mapping $m \mapsto \Vert m \Vert = \int_{G} \chi _{G}d \vert m \vert$ is a norm on $M_{1}(G,A)$, where $\chi _{G}$ denotes the characteristic function of $G$.
\end{definition}
In the sequel, $K$ will designate a compact subgroup of $G$, $\nu$ and $\lambda$  left Haar measures respectively on $K$ and $G$.
\subsection{G-invariant measure}
\begin{definition}\mbox{ }\\
Let $\mu$ be  a Radon measure on $G/K,$ the homogenous space of left $K$-cosets and $g$  an element of $G$. Define $\mu_{g}$ by $\mu_{g}(E) = \mu(gE)$ for Borel subsets $E$ of $G/K$. The measure $\mu$ is called $G-$invariant measure if $\mu_{g} = \mu$, for $g \in G$. See \cite{folland}, \cite{gaal} for more details.
\end{definition}
 
Throughout the paper, $ \mu $  designates the $ G- $invariant  measure on $G/K$.
\begin{theorem}\mbox{ }\\
 For any $ f \in \mathcal {K} (G), $ we have  \cite{folland}, \cite{anahar}:
\begin{eqnarray} \label{inv}
\int_{G}f(g)d\lambda(g) &=& \int_{X}d\mu(\dot{g})\int_{K}f(gk)d\nu(k).
\end{eqnarray} 
The previous formula (\ref{inv}) extends also to  every $f\in L_{1}(G,\lambda,A)$ \cite{Bar2}.\\
\end{theorem}
\subsection{Group representation}
\begin{definition}\mbox{ }\\
 A unitary representation of $G$ is  a homomorphism $L$ from $G$ into the group $U (H)$ of the unitary linear operators on some nonzero  Hilbert space $H,$ which is continuous with respect to
the strong operator topology satisfying: for   $g_{1}, g_{2} \in G$,
\[L_{( g_{1}, g_{2})}  =  L_{g_{1}}L_{g_{2}} \quad \text{and  } \quad L_{e}  =  Id_{H}.\]
$H$ is called the
representation space of $L$, and its dimension is called the dimension or degree of $L$.\\
Suppose $\mathcal{M}$ is a closed subspace of $H$. $\mathcal{M}$ is called an invariant subspace for $L$ if $L_{g}\mathcal{M} \subset \mathcal{M}$ $\forall g \in G$. If $\mathcal{M}$ is invariant and $\mathcal{M} \neq \{ 0\},$ then $L^{\mathcal{M}}$ such that
 \[L^{\mathcal{M}}_{g} = L_{g} \vert_{\mathcal{M}}\]
 defines a representation of $G$ on $\mathcal{M}$, called a subrepresentation of $L$. 
 If $L$ admits an invariant subspace that is nontrivial (i.e. $\neq \{0\}$ or $H$) then $L$ is called reducible, otherwise $L$ is irreducible. If $G$ is compact and $L$ irreducible then the dimension of $L$ is finite.
\end{definition} 
\begin{definition} \mbox{ }\\
Two unitary irreductible representations $L$ and $V$ into
 $H$ and $N,$ respectively, are said to be equivalent if there is an isomorphism
$T :  H\longrightarrow N $ such that, $\forall t \in G,$
\[  T\circ L_{t} = V_{t}\circ T.\]
\end{definition}
Consider now the subgroup $K$, $\Sigma$ the coset space (called the dual object of $ K $), $\sigma\in \Sigma,$ $L^{\sigma}$ a representative of $\sigma$, $H_{\sigma}$ a representation space of $L^{\sigma},$ and $d_{\sigma}$ its dimension.
\begin{theorem}\mbox{ }\\
Let  $(L_{ij}^{\sigma})_{1 \leq i,j \leq d_{\sigma}}$ be the matrix of $L^{\sigma}$ in an orthonormal  basis  $(\xi_{i})_{i=1}^{d\sigma}$ of $H_{\sigma}.$ Then,  (see \cite{anahar}, \cite{men}, \cite{folland}, \cite{B-R}), 
 \begin{eqnarray} \label{sch1}
 \int_{K}L_{ij}^{\sigma}(t)\overline{L_{lm}^{\sigma}}d \nu (t) & = & \dfrac{\delta _{il} \delta_{jm}}{d_{\sigma} }
\end{eqnarray}  
and
 \begin{eqnarray} \label{sch2}
 \int_{K}L_{ij}^{\sigma}(t)\overline{L_{lm}^{\tau}}d \nu (t) & = & 0 \text{ if } \sigma \neq \tau.
\end{eqnarray}
\end{theorem}  
\subsection{Induced representation} \mbox{ }\\
Let  $q : G \longrightarrow G/K$ be the canonical quotient map of $G$ into $G/K$ and suppose $H_{\sigma}$ separable.
Denote by  $H^{L^{\sigma}}_{ 0}$ the set
 \begin{eqnarray}
H^{L^{\sigma}}_{ 0} = \left\lbrace u\in \mathfrak{C}(G,H_{\sigma})\mbox{ } : \mbox{ } q(\mbox{Supp(u))}\text{ is compact } \text{ and  } u(gk) = L^{\sigma}_{k^{-1}}u(g)\right\rbrace.
\end{eqnarray}
\begin{proposition}\mbox{ }\\
If $\eta :G \longrightarrow H_{\sigma}$ 
is  continuous  with compact support, then the function $u_{\eta}$ such that
\begin{eqnarray} \label{form}
u_{\eta}(g) = \int_{K}L^{\sigma}_{k}\eta(gk)d\nu(k)
\end{eqnarray} 
belongs to $H^{L^{\sigma}}_{ 0},$ and is uniformly continuous on $G$. Moreover, every 
element of $H^{L^{\sigma}}_{ 0}$ is of the form $u_{\eta}$. See \cite{folland} for more details.
\end{proposition}
\begin{proposition}\mbox{ }\\
The mapping: 
\begin{eqnarray}
(u,v)&\longmapsto & \left<u,v \right> = \int_{G/K}\left<u(g),v(g)\right>_{H_{\sigma}}d\mu(\dot{g}).
\end{eqnarray}
on  $H^{L^{\sigma}}_{0} \times H^{L^{\sigma}}_{0}$ is an inner product on $H^{L^{\sigma}}_{0}.$ 
\end{proposition}
\begin{proposition}\mbox{ }\\
$ G$ acts on $H^{L^{\sigma}}_{0}$ by left translation, $u \longmapsto L_{t}u$, so we obtain a unitary
representation of $G$ with respect to this inner product on $H^{L^{\sigma}}_{0}$.
The inner product is preserved by left translations, since $\mu$ is invariant. Hence, if we denote by $H^{L^{\sigma}}$ the Hilbert space completion of $H^{L^{\sigma}}_{0},$ the translation operators $L_{t}$  extend to unitary operators on $H^{L^{\sigma}}$. Then the map $t \longmapsto L_{t}u$ is continuous from $G$ to $H^{L^{\sigma}}$ for each $u \in H^{L^{\sigma}}_{0},$ and then, since the operators $L_{t}$ are uniformly bounded, they are strongly continuous on $H^{L^{\sigma}}.$ Hence they define a unitary representation of G, called the representation induced by $L^{\sigma},$  denoted by $U^{L^{\sigma}}:$
\[ U_{t}^{L^{\sigma}}u(g) = L_{t}u(g) = u(t^{-1}g).\]
The representation space is $H^{L^{\sigma}}$.
\end{proposition}
\begin{remark}\mbox{ }\\
The representations of $G$ induced from $K$ are generally infinite-dimensional unless $G/K$ is a finite set.\\
The induced representation $U^{L^{\sigma}}$ is irreducible only if $L^{\sigma}$ is irreducible \cite{gaal}.
\end{remark} 
\subsection{Fourier-Stieltjes transform defined by induced representation on locally compact groups } \cite{yao}

We have $K$ compact, $H_{\sigma}$ is finite and separable;  then $H^{L^{\sigma}}$ is also separable.  Each of these spaces admits a Hilbertian basis according to Gram-Schmidt process.\\\\
In this work, we suppose that $K$ is chosen such that $L^{\sigma}$ and $U^{L^{\sigma}}$ are both irreducible.
\\
Let $(\theta_{i})_{i=1}^{\infty}$ be an orthonormal basis of $H^{L^{\sigma}}$ and  $ (\xi_{i})_{i=1}^{d\sigma}$ be an orthonormal  basis of $H_{\sigma}$.
We define
\begin{eqnarray} \label{trig}
u_{ij}^{L^{\sigma}}(t) &:=& \left< U_{t}^{L^{\sigma}}\theta_{j},\theta_{i}\right>_{H^{L^{\sigma}}}\cr
 &=& \int_{G/K}\left< \theta_{j}(t^{-1}g), \theta_{i}(g)\right>d\mu(\dot{g})
\end{eqnarray}
and 
\[L_{ij}^{\sigma}(k) := \left< L_{k}^{\sigma}\xi_{j},\xi_{i}\right>_{H_{\sigma}}
\]
As a result, there is a family of mappings $ (\alpha_{is})_{s=1}^{d\sigma}$ of $G$ into $\mathbb{K} with $ ($\mathbb{K} = \mathbb{R}$ or $ \mathbb{C}$) such that \[\theta _{i}(g) = \sum _{s = 1}^{d \sigma} \alpha_{is}(g)\xi _{s}.\]
\begin{theorem} \mbox{   }
We have
\begin{eqnarray}
\int_{G}u_{ij}^{L^{\sigma}}(t)\overline{u}_{lm}^{L^{\sigma}}(t)d\lambda(t) & = & \dfrac{c_{ijlm}}{d_{\sigma}}
\end{eqnarray}
where  \begin{eqnarray*}
 c_{ijlm} & : = &   \int _{G/K} d\mu(\dot{t})\sum_{r,s = 1}^{d_{\sigma}} \int_{G/K}\alpha_{js}(t^{-1}g)\overline{\alpha}_{ir}(g)d\mu(\dot{g})\int_{G/K}\alpha_{ms}(t^{-1}h)\overline{\alpha}_{lr}(h)d\mu(\dot{h})\\
 &  = &    \sum_{r,s = 1}^{d_{\sigma}} \int_{(G/K)^{3}} \alpha_{js}(t^{-1}g)\overline{\alpha}_{ir}(g) \alpha_{ms}(t^{-1}h)\overline{\alpha}_{lr}(h)d\mu(\dot{g})d\mu(\dot{h})d\mu(\dot{t})
 \end{eqnarray*}
and
\begin{eqnarray}
\int_{G}u_{ij}^{L^{\sigma}}(t)\overline{u}_{lm}^{L^{\tau}}(t)d\lambda(t) & = & 0 \text{ if } \sigma \neq \tau.
\end{eqnarray}
\end{theorem}
In the case of a particular orthonormal basis, the orthogonality relation reduces to the following:
\begin{corollary} \mbox{ }
Choosing  an orthonormal  basis $(\xi_{i})_{i=1}^{d\sigma}$  of $H_{\sigma}$ such that 
\begin{eqnarray}\label{normal}
\int_{(G/K)^{3}}\alpha_{js}(t^{-1}g)\overline{\alpha}_{ir}(g) \alpha_{ms}(t^{-1}h)\overline{\alpha}_{lr}(h)d\mu(\dot{g})d\mu(\dot{h})d\mu(\dot{t}) = \left\{\begin{array}{cl}
 \dfrac{1}{d_{\sigma}^{2}}  &\text{  if } j = m \text{  and } i = l \\
  & \\
0 &\text{ if not}
\end{array}\right. 
\end{eqnarray}
(or, also, $c_{ijlm} = \delta_{il}\delta_{jm}$)
leads to 
\begin{eqnarray}
\int_{G}u_{ij}^{L^{\sigma}}(t)\overline{u}_{lm}^{L^{\sigma}}(t)d\lambda(t) & = & \dfrac{\delta_{il}\delta_{jm}}{d_{\sigma} }.
\end{eqnarray}
\end{corollary}

In the sequel, $(\theta_{i})_{i=1}^{\infty}$ will designate an orthonormal basis of $H^{L^{\sigma}}$ and  $ (\xi_{i})_{i=1}^{d\sigma}$ that of $H_{\sigma},$ where $ (\xi_{i})_{i=1}^{d\sigma}$ is chosen such that
\begin{eqnarray*}
\int_{G/K}d\mu(\dot{t})\int_{G/K}\alpha_{js}(t^{-1}g)\overline{\alpha}_{ir}(g)d\mu(\dot{g})\int_{G/K}\alpha_{ms}(t^{-1}h)\overline{\alpha}_{lr}(h)d\mu(\dot{h}) = \left\{\begin{array}{cl}
 \dfrac{1}{d_{\sigma}^{2}}  &\text{  if } j = m \text{  and } i = l \\
  & \\
0 &\text{ if not.}
\end{array}\right. 
\end{eqnarray*}
\begin{definition}\label{deffour1} 
Assume $m\in M_{1}(G,\mathcal{A})$.
The Fourier-Stieltjes transform of an arbitrary measure $m$ is defined as the family $(\hat{m}(\sigma))_{\sigma\in \Sigma}$ of sesquilinear mappings on $H^{L^{\sigma}} \times H^{L^{\sigma}}$ into  $\mathcal{A}$, given by the relation
\begin{eqnarray}
\hat{m}(\sigma)(u,v) =  \int_{G}\left< \overline{U}_{t}^{L^{\sigma}}u,v \right>_{H^{L^{\sigma}}}dm(t) = \int_{G}\int_{G/K} \left<\overline{u}(t^{-1}g),v(g)\right>d\mu(\dot{g})dm(t).
\end{eqnarray}
\end{definition}
\begin{definition}\label{deffour2}
 For any $f\in L_{1}(G,\mathcal{A},\lambda)$, where $\lambda$ denotes a Haar measure on $ G $, the Fourier transform of $f$ is a family $(\hat{f}(\sigma))_{\sigma\in \Sigma}$ of sesquilinear mappings of $H^{L^{\sigma}} \times H^{L^{\sigma}}$  into $\mathcal{A},$ given by the relation
\begin{eqnarray}
\hat{f}(\sigma)(u,v) =  \int_{G}\left< \overline{U}_{t}^{L^{\sigma}}u,v \right>_{H^{L^{\sigma}}}f(t)d\lambda(t) = \int_{G}\int_{G/K} \left<\overline{u}(t^{-1}g),v(g)\right>f(t)d\mu(\dot{g})d\lambda(t).
\end{eqnarray}
\end{definition}
Let denote by $\mathcal{S}(\Sigma,\mathcal{A}) = \displaystyle{\prod_{\sigma\in\Sigma}}\mathcal{S}(H^{L^{\sigma}}\times H^{L^{\sigma}}, \mathcal{A}),$ where $\mathcal{S}(H^{L^{\sigma}}\times H^{L^{\sigma}}, \mathcal{A})$ is the set of all the sesquilinear mappings of $H^{L^{\sigma}}\times H^{L^{\sigma}}$ into $\mathcal{A}$.

Also, assume $ \displaystyle{\prod_{\sigma\in\Sigma}}\mathcal{S}(H^{L^{\sigma}}\times H^{L^{\sigma}}, \mathcal{A})$ is a vector space for addition and multiplication by a scalar of mappings. 

\begin{proposition}\label{dec} \mbox{  } 
Let us consider $\Phi$ in  $\mathcal{S}(\Sigma,\mathcal{A}).$ 
There is the matrix $(a_{ij}^{\sigma})_{1 \leq i, j \leq \infty}$,  $a_{ij}^{\sigma}\in \mathcal{A}$ such that 
\[\Phi (\sigma) = \sum _{i,j = 1}^{\infty}d_{\sigma}  a_{ij}^{\sigma}\hat{u}_{ij}^{L^{\sigma}} \] 
with $\hat{u}_{ij}^{L^{\sigma}}$ the Fourier transform of $u_{ij}^{L^{\sigma}}$ given by
\[\hat{u}_{ij}^{L^{\sigma}}(\sigma)(u,v)  =  \int_{G}\left< \overline{U}_{t}^{L^{\sigma}}u,v \right>_{H^{L^{\sigma}}}u_{ij}^{L^{\sigma}}(t)d\lambda(t).\]
\end{proposition}
\begin{corollary}
For any $f \in L_{1}(G, \mathcal{A}, \lambda)$, there is a matrix  $(a_{ij})_{1 \leq i, j \leq \infty}$, $a_{ij}\in \mathcal{A}$  such that  
\[\hat{f} (\sigma) = \sum _{i,j = 1}^{\infty}d_{\sigma} a_{ij}^{\sigma}\hat{u}_{ij}^{L^{\sigma}}.  \]
\end{corollary}

\begin{theorem}\label{inj} The mapping $ m \longmapsto \hat{m}$ from $ M_{1}(G,\mathcal{A})$ into $\mathcal{S}_{\infty}(\Sigma,\mathcal{A})$ is linear,  injective and continuous, where $$\mathcal{S}_{\infty}(\Sigma,\mathcal{A})=\left\lbrace \Phi \in \mathcal{S}(\Sigma,\mathcal{A}): \Vert \Phi \Vert_{\infty} < \infty 
\right\rbrace.$$
\end{theorem}
\begin{corollary}\label{inj2}
The mapping $f \longmapsto \hat{f}$ from $ L_{1}(G, \lambda, \mathcal{A})$ into $\mathfrak{S}_{\infty}(\Sigma,\mathcal{A})$ is linear, injective, and continuous.
\end{corollary}  

\begin{theorem}
For every $f\in L_{2}(G,\mathcal{A})$, there is $a_{ij}^{\sigma} \in \mathcal{A},$  $ 1 \leq i, j < \infty$, $\sigma \in \Sigma$ such that 
\begin{eqnarray}
f & = & \displaystyle\sum_{\sigma \in \Sigma}d_{\sigma}\displaystyle\sum_{i,j = 1}^{\infty}a_{ij}^{\sigma}u_{ij}^{L^{\sigma}}.
\end{eqnarray} 
\end{theorem}
\begin{corollary}
For every $f \in L_{2}(G,\mathcal{A})$, we have 
\begin{eqnarray}
f  & = & \displaystyle\sum_{\sigma \in \Sigma}d_{\sigma}\displaystyle\sum_{i,j = 1}^{\infty}\hat{f}(\theta_{j},\theta_{i}) u_{ij}^{L^{\sigma}}. 
\end{eqnarray}
\end{corollary}
\subsection{Tensor Products of Vector Spaces} \cite{ryan} \mbox{ }\\
We recall that a mapping $D$ from the cartesian product $E \times F$ of vector
spaces into the vector space $S$ is bilinear if it is linear in each variable, that
is,
\begin{eqnarray*}
D(\alpha_{1}x_{1} +\alpha_{2}x_{2},y) & = & \alpha_{1}D(x_{1},y) + \alpha_{2}D(x_{2},y) \text{ and }\\
D(x,\beta_{1}y_{1} +\beta_{2}y_{2}) & = & \beta_{1}D(x, y_{1}) + \beta_{2}D(x, y_{2})
\end{eqnarray*}
for all $x_{i},x \in E$, $y_{i}, y\in F$ and all scalars $\alpha_{i}, \beta_{i}$, $i = 1, 2$.\\
We write $B(E \times F, S)$ for the vector space of bilinear mappings from the product $E \times F$ into $S$; when $S$ is the scalar field we denote the corresponding space of bilinear forms simply by $B(E \times F)$.
\begin{definition} \mbox{ }\\
The tensor product, $E \otimes F$, of the vector spaces $E$, $F$ can be constructed as a space of linear functionals on $B(E \times F)$, in the following way: for $x\in E$, $y\in F$, we denote by $x \otimes y$ the functional given by evaluation at the point $(x,y)$. In other words,
\[ (x \otimes y) (D) = \left < D, x \otimes y \right > = D(x,y) \] for each bilinear form $D$ of $E \times F$. The tensor product $E \otimes F$ is the subspace of the dual $B(E \times F)^{*}$ spanned by these elements. Thus, a typical tensor in $E \otimes F$ has the form
\begin{eqnarray}
u & = & \Sigma_{i = 1}^{n}\lambda_{i} x_{i} \otimes y_{i}
\end{eqnarray}
where $n$ is a natural number, $\lambda_{i} \in \mathbb{K}$, $x_{i}\in E$ and $y_{i} \in F$.
\end{definition}
\begin{remark} \mbox{ }\\
It is important to realize that the representation of $u$ is not unique - in general, there will be many different ways to write a given tensor in the above form.
\end{remark}
\begin{proposition} \mbox{ }\\
Let $E$ and $F$ be Banach spaces. The mapping $\pi$ defined by
\begin{eqnarray}
\pi (u)  & = & \inf \{ \Sigma_{i = 1}^{n} \Vert x_{i} \Vert_{E}\Vert y_{i}\Vert _{F}; \mbox{  } u = \Sigma_{i = 1}^{n} x_{i}\otimes y_{i} \}
\end{eqnarray}
is a norm in $E\otimes F$ known as the projective norm.
\end{proposition}
We shall denote by $E \otimes_{\pi} F$ tensor product $E \otimes F$ endowed with
the projective norm, $\pi$. Unless the spaces $E$ and $F$ are finite dimensional,
this space is not complete. We denote its completion by $E \hat{\otimes}_{\pi} F$.
The Banach space $E \hat{\otimes}_{\pi} F$ will be referred to as the projective tensor product of the Banach of $E$ and $F$.
\begin{proposition} \mbox{ }\\
 Let $S$ and $A$ be  Banach spaces and $\mu$ any measure on $S$. $L_{1}(S,\mu, \mathbb{K})\hat{\otimes}_{\pi} =   L_{1}(S,\mu, A)$.
\end{proposition}
\begin{proposition} \mbox{ }\\
Let $E$ and $F$ be Banach spaces. The mapping $\varepsilon$ defined by
\begin{eqnarray}
\varepsilon (u)  & = & \sup \{ \vert \Sigma_{i = 1}^{n} \varphi ( x_{i}) \psi (y_{i})\vert; \mbox{  } \varphi \in E^{*}, \psi \in F^{*}, \Vert \varphi \Vert \leq 1, \Vert \psi \Vert, u = \Sigma_{i = 1}^{n} x_{i}\otimes y_{i} \}
\end{eqnarray}
is a norm in $E\otimes F$ known as the injective norm.
\end{proposition}
We shall denote by $E \otimes_{\varepsilon} F$ tensor product $E \otimes F$ endowed with
the injective norm, $\pi$. Unless the spaces $E$ and $F$ are finite dimensional,
this space is not complete. We denote its completion by $E \hat{\otimes}_{\varepsilon} F$.
The Banach space $E \hat{\otimes}_{\varepsilon} F$ will be referred to as the injective tensor product of the Banach of $E$ and $F$.
\begin{proposition}\label{tens1} \mbox{ }\\
 Let $S$ and $A$ be  Banach spaces. $\mathfrak{C}(S, \mathbb{K})\hat{\otimes}_{\varepsilon} =   \mathfrak{C}(S, A)$.
\end{proposition}
\subsection{Stone-Weierstrass theorem for locally compact spaces} \mbox{ }\\
\begin{theorem} \label{stone} \mbox{  }\\
Let $\Omega$ be a locally compact space, which is non-compact, let $ \mathfrak{C}_{0}^{\mathbb{K}}(\Omega)$ be the algebra of continuous $\mathbb{K}$-valued functions on $\Omega$ that vanish at infinity, as above, equipped with the supremum norm and let $\mathcal{V} \subset \mathfrak{C}_{0}^{\mathbb{K}}(\Omega)$ with the following separation properties:
\begin{itemize}
\item[(i)] for any two points $\omega_{1}, \omega_{2} \in \Omega$ with $\omega_{1}\neq \omega_{2}$  there exists $f \in\mathcal{V}$ such that $f(\omega_{1}) \neq f(\omega_{2})$
\item[(ii)] for any $\omega \in \Omega$ there exists $f \in\mathcal{V}$ with $f(\omega) \neq0$.\\
\item[(a)] If $\mathbb{K} = \mathbb{R}$ then  $\in\mathcal{V}$ is dense in $\mathfrak{C}_{0}^{\mathbb{R}}(\Omega).$ 
\item[(b)] If $\mathbb{K} = \mathbb{C}$ and if $\mathcal{V}$ is a $_{*}$-subalgebra (i.e.  $f\in \mathcal{V}\Longrightarrow \overline{f}\in \mathcal{V}$) then $\mathcal{V}$ is dense in $\mathfrak{C}_{0}(\Omega).$  
\end{itemize}
\end{theorem}

  
\section{Main results}
\subsection{$\mathfrak{S}_{p} $  spaces}
In this section we define some particular spaces of Fourier transform and give some properties of them.

\begin{definition}
\begin{itemize}
\item[(i)] $\displaystyle \mathfrak{S}_{\infty}(\Sigma,\mathcal{A})$ =  $\left\lbrace \Phi \in \mathfrak{S}(\Sigma,\mathcal{A}): \Vert \Phi \Vert_{\infty} < \infty 
\right\rbrace$.  

\item[(ii)] $\displaystyle \mathfrak{S}_{00}(\Sigma,\mathcal{A})$ = $\left\lbrace \Phi \in \mathfrak{S}(\Sigma,\mathcal{A}): \left\lbrace  \sigma \in\Sigma :\Phi(\sigma) \neq 0 \right\rbrace \quad \mbox{is finite}\right\rbrace$.
\item[(iii)] $\displaystyle\mathfrak{S}_{0}(\Sigma,\mathcal{A})$ = $\left\lbrace \Phi \in \mathfrak{S}(\Sigma,\mathcal{A}): \forall \varepsilon > 0 \quad\lbrace  \sigma \in\Sigma :\quad \Vert\Phi(\sigma)\Vert > \varepsilon \rbrace \quad\mbox{is finite}\right\rbrace$.
\item[(iv)] $ \displaystyle\mathfrak{S}_{p}(\Sigma,\mathcal{A})  =  \left\lbrace \Phi \in \mathfrak{S}(\Sigma,\mathcal{A}): \quad \sum_{\sigma\in\Sigma}d\sigma\sum_{i,j = 1}^{\infty}\Vert\Phi(\sigma)(\theta_{i}^{\sigma},\theta_{j}^{\sigma}) \Vert_{\mathcal{A}}^{P} < \infty\right\rbrace$. 
\end{itemize} 
\end{definition}
\begin{proposition}
$$\mathfrak{S}_{00}(\Sigma,\mathcal{A}) \subset \mathfrak{S}_{1}(\Sigma,\mathcal{A})\subset...\subset \mathfrak{S}_{p}(\Sigma,\mathcal{A})\subset ... \subset \mathfrak{S}_{0}(\Sigma,\mathcal{A})\subset \mathfrak{S}_{\infty}(\Sigma,\mathcal{A}).$$ 
\end{proposition} 
\textbf{Proof} \\
\begin{itemize}
\item[(i)] Let us first show that for $1 \leq p \leq q < \infty$, $\mathfrak{S}_{p}(\Sigma,\mathcal{A}) \subset \mathfrak{S}_{q}(\Sigma,\mathcal{A})$.
For this, we  show that $\Vert . \Vert_{q} \leq \Vert . \Vert _{p},$ where 
\[\Vert \Phi \Vert_{\varphi} = \left(\sum_{\sigma \in \Sigma}d\sigma \sum_{i,j=1}^{\infty}\left\Vert \Phi (\sigma)(\theta_{i}^{\sigma},\theta_{j}^{\sigma}) \right\Vert_{\mathcal{A}}^{\varphi}\right)^{\frac{1}{\varphi}}.   \]
Put \[\varrho = \Vert \Phi \Vert_{q} = \left(\sum_{\sigma \in \Sigma}d\sigma \sum_{i,j=1}^{\infty}\left\Vert \Phi (\sigma)(\theta_{i}^{\sigma},\theta_{j}^{\sigma}) \right\Vert_{\mathcal{A}}^{q}\right)^{\frac{1}{q} } \]
and \[b_{\sigma} =  \sum_{i,j = 1}^{\infty}\Vert\Phi(\sigma)(\theta_{i}^{\sigma},\theta_{j}^{\sigma}) \Vert_{\mathcal{A}}.\]
The  summation in the  $\varrho ^{q}$ expression can permit to  retrieve the term $d_{\sigma}b_{\sigma}^{q}$.
Then, we have  $1 \leq d_{\sigma}b_{\sigma}^{q} \leq \varrho^{q} \Longrightarrow 1 \leq d_{\sigma}^{\frac{1}{q}}d_{\sigma} \leq \varrho$.
Thus, 
\[0 < \dfrac{d_{\sigma}^{\frac{1}{q}}d_{\sigma}}{\varrho} \leq 1, \]
and we have:
\begin{eqnarray*}
\left(\dfrac{d_{\sigma}^{\frac{1}{q}}d_{\sigma}}{\varrho}\right)^{q} \leq \left(\dfrac{d_{\sigma}^{\frac{1}{q}}d_{\sigma}}{\varrho}\right)^{p} &\Longrightarrow & \sum_{\sigma \in \Sigma}\left(\dfrac{d_{\sigma}^{\frac{1}{q}}d_{\sigma}}{\varrho}\right)^{q} \leq \sum_{\sigma\in\Sigma}\left(\dfrac{d_{\sigma}^{\frac{1}{q}}d_{\sigma}}{\varrho}\right)^{p}\\
& \Longrightarrow & \sum_{\sigma \in \Sigma}\left(\dfrac{d_{\sigma}b_{\sigma}^{q}}{\varrho^{q}}\right) \leq \displaystyle{\sum_{\sigma\i\Sigma}\left(\dfrac{d_{\sigma}^{\frac{p}{q}}b_{\sigma}^{p}}{\varrho^{p}}\right)}\\
& \Longrightarrow &  \dfrac{\displaystyle{\sum_{\sigma\in\Sigma}d_{\sigma}b_{\sigma}^{q}}}{\varrho^{q}} \leq \dfrac{\displaystyle{\sum_{\sigma\in\Sigma} d_{\sigma}^{\frac{p}{q}}b_{\sigma}^{p}}}{\varrho^{p}}\\
& \Longrightarrow & 1 \leq \dfrac{\displaystyle{\sum_{\sigma\in\Sigma} d_{\sigma}^{\frac{p}{q}}b_{\sigma}^{p}}}{\displaystyle{\sum_{\sigma\in\Sigma}d_{\sigma}b_{\sigma}^{q}}}\\
& \Longrightarrow & \sum_{\sigma\in\Sigma}d_{\sigma}b_{\sigma}^{q} \leq \sum_{\sigma\in\Sigma}d_{\sigma}^{\frac{p}{q}}b_{\sigma}^{p}
\end{eqnarray*}
as $d_{\sigma}\geq 1$ and $\frac{p}{q} \leq 1$. Then, $d_{\sigma}^{\frac{p}{q}} \leq d_{\sigma},$ and
\begin{eqnarray*}
\sum_{\sigma\in\Sigma}d_{\sigma}b_{\sigma}^{q} \leq \sum_{\sigma\in\Sigma}d_{\sigma}^{\frac{p}{q}}b_{\sigma}^{p} \leq \sum_{\sigma\in\Sigma}d_{\sigma}b_{\sigma}^{p}.
\end{eqnarray*}
Hence $\Vert \Phi\Vert_{q} \leq \Vert \Phi \Vert_{p}$.\\ \\
\item[(ii)] For $q = \infty$, $p < \infty$ and $\Phi  \in  \mathfrak{S}_{p}(\Sigma,\mathcal{A}),$ \\ 
\begin{eqnarray*}
\Vert \Phi \Vert_{\infty} & = & \sup \lbrace \Vert \Phi (\sigma) \Vert / \sigma \in\Sigma \rbrace
 \leq
 \left(\sum_{\sigma \in \Sigma} d_{\sigma}\sum_{i,j = 1}^{\infty} \Vert \Phi(\sigma)(\theta_{j},\theta_{i}) \Vert_{\mathcal{A}}^{p} \right)^{\frac{1}{p}}. 
\end{eqnarray*} 
\item[(iii)]  One  shows that for $\Phi \in \mathfrak{S}_{0}(\Sigma,\mathcal{A}), $ $ \Vert \Phi \Vert_{\infty} < \infty$ implying that $ \Phi$ belongs to $ \mathfrak{S}_{\infty}(\Sigma,\mathcal{A})$. Then we have  $ \mathfrak{S}_{0}(\Sigma,\mathcal{A})\subset  \mathfrak{S}_{\infty}(\Sigma,\mathcal{A}).$ \\
\item[(iv)] For any $p$, $1 \leq p < \infty$ and $\Phi  \in  \mathfrak{S}_{p}(\Sigma,\mathcal{A}),$ we have $\displaystyle{\sum_{\sigma\in\Sigma}d\sigma\sum_{i,j = 1}^{\infty}\Vert\Phi(\sigma)(\theta_{i}^{\sigma},\theta_{j}^{\sigma}) \Vert_{\mathcal{A}}^{P} < \infty }$. It means that $\forall \varepsilon > 0,$ the set $\lbrace  \sigma \in\Sigma  \mbox{ }\mbox{ }:\mbox{ }\mbox{ }\Vert\Phi(\sigma)\Vert > \varepsilon \rbrace \mbox{ }\mbox{ }\mbox{ is finite}.$ Then,   $\mathfrak{S}_{p}(\Sigma,\mathcal{A}) \subset   \mathfrak{S}_{0}(\Sigma,\mathcal{A})$. \\
\item[(v)] For any $ \Phi \in  \mathfrak{S}_{00}(\Sigma,\mathcal{A}),$ we have $\displaystyle{\sum_{\sigma\in\Sigma}d\sigma\sum_{i,j = 1}^{\infty}\Vert\Phi(\sigma)(\theta_{i}^{\sigma},\theta_{j}^{\sigma}) \Vert_{\mathcal{A}}^{P} < \infty }$ for any $p$, $1 \leq p < \infty$ because the set  $\lbrace  \sigma \in\Sigma  \mbox{ }\mbox{ }/\mbox{ }\mbox{ }\Phi(\sigma) \neq 0 \rbrace \mbox{ }\mbox{ }\mbox{ is finite}.$ Then,  $\mathfrak{S}_{00}(\Sigma,\mathcal{A}) \subset   \mathfrak{S}_{p}(\Sigma,\mathcal{A}). \blacktriangledown$
\end{itemize} 
\begin{proposition}
 \mbox{ }
$\mathfrak{S}_{00}(\Sigma,\mathcal{A})$ is  dense in  $\mathfrak{S}_{0}(\Sigma,\mathcal{A})$.
\end{proposition}
\textbf{Proof} 
Let $\phi$ be an element of $\mathfrak{S}_{0}(\Sigma,\mathcal{A})$. Consider the sequence $( \phi _{n})_{n\in \mathbb{N}}$ defined by \\
$ \phi _{n} (\sigma)= $
$\left \{
\begin{array}{rcl}
\phi(\sigma)  \mbox{ if  }  \Vert  \phi(\sigma) \Vert \geq \dfrac{1}{n} \\ \mbox{   }  \\ 0  \mbox{     if   } \Vert  \phi(\sigma) \Vert < \dfrac{1}{n}. 
\end{array}
\right. $\\
We can show that for any $n \in \mathbb{N}$, $\phi _{n}$ belongs to  $ \mathfrak{S}_{00}(\Sigma,\mathcal{A}),$ because, if  $ \phi _{n} (\sigma) \neq 0,$ it means that $ \phi _{n} (\sigma) = \phi  (\sigma),$  and then $ \Vert  \phi(\sigma) \Vert \geq \dfrac{1}{n}$. Since $ \phi$ belongs  to  $\mathfrak{S}_{0}(\Sigma,\mathcal{A})$, the set $\lbrace \sigma \in \Sigma \mbox{      }  :  \Vert  \phi(\sigma) \Vert \geq \dfrac{1}{n} \rbrace$  is finite,  implying the set $\lbrace \sigma \in \Sigma \mbox{      }  :   \phi_{n}(\sigma)  \neq 0 \rbrace$  is finite.\\ \\
$ \Vert \phi_{n} - \phi \Vert _{\infty} \leq \Vert \phi \Vert _{\infty} < \dfrac{1}{n}$.\\
For any $\varepsilon > 0$ there exists $N$ such that $\dfrac{1}{N} < \varepsilon. $
Then for any $n > N,$ we have $\dfrac{1}{n} < \dfrac{1}{N} < \varepsilon$.\\
Thus $ \Vert \phi_{n} - \phi \Vert _{\infty} \leq \Vert \phi \Vert _{\infty} < \dfrac{1}{n} < \dfrac{1}{N} < \varepsilon.$ $\blacktriangledown$
\subsection{Lesbesgue Theorem}
Throughout this section, we denote by $ \hat{X}$ the set $\{ \hat{f}: f \in X \}.$
For any $\sigma \in \Sigma,$ let us denote by $\mathcal{L}_{\sigma}(G)$ the set of linear combinations of functions $t  \longmapsto <\overline{U_{t}^{L^{\sigma}}}\theta_{j};\theta_{i}>_{H^{L^{\sigma}}},$ and by  $\mathcal{L}(G) = \displaystyle\bigcup_{\sigma \in \Sigma } \mathcal{L}_{\sigma}(G).$
\begin{lemma} \mbox{  }\\
\begin{itemize}
\item[(i)] $\widehat{\mathcal{L}_{\sigma}(G)\otimes \mathcal{A}} = S(H^{L^{\sigma}}\times H^{L^{\sigma}}, \mathcal{A})$  and 
\item[(ii)] $ \widehat{\mathcal{L}(G)\otimes \mathcal{A}} = \mathfrak{S}_{00}(\Sigma,\mathcal{A}).$
\end{itemize}
\end{lemma}
\textbf{Proof}\\
\begin{itemize}
\item[(i)] According to  Propostion \ref{dec}, $\Phi(\sigma) \in S(H^{L^{\sigma}}\times H^{L^{\sigma}}, A) \Longleftrightarrow$  there exists $a_{ij} \in A$ and $u_{ij}^{L^{\sigma}} \in \mathcal{L}_{\sigma}(G)$ such that  $\Phi (\sigma) = \displaystyle\sum _{i,j = 1}^{\infty}d_{\sigma}  a_{ij}^{\sigma}\hat{u}_{ij}^{L^{\sigma}}$.
According to theorem \ref{inj} we have $\Vert \Phi \Vert_{\infty} < \infty$ then there exists $n\in\mathbb{N}$ such that $\displaystyle\sum _{i,j = 1}^{n}d_{\sigma}  a_{ij}^{\sigma}\hat{u}_{ij}^{L^{\sigma}} \neq 0.$ Thus we have $\Phi (\sigma) = \displaystyle\sum _{i,j = 1}^{n}d_{\sigma}  a_{ij}^{\sigma}\hat{u}_{ij}^{L^{\sigma}}$.\\ Then $\Phi(\sigma)\in \widehat{\mathcal{L}_{\sigma}(G)\otimes A}.$
\item[(ii)] Now consider $f \in\mathcal{L}(G)\otimes A,$ there exists $n\in \mathbb{N}$ such  that $f = \displaystyle \sum _{k = 1}^{n} \mu_{k} f_{\sigma_{k}}$, where $\mu_{k}\in \mathbb{C},\sigma_{k} \in \Sigma$ and $f_{\sigma_{k}} = \displaystyle \sum_{i,j}x_{ij}u_{ij}^{L^{\sigma_{k}}}$ with $x_{ij} \in A$ and $u_{ij}^{L^{\sigma_{k}}} \in \mathcal{L}_{\sigma_{k}}(G).$ \\
$\hat{f}(\sigma)(u,v) = \displaystyle \sum _{k = 1}^{n} \mu_{i}\displaystyle \sum_{i,j}x_{ij}\hat{u}_{ij}^{L^{\sigma_{k}}}(\sigma)(u,v)$ then $\hat{f}(\sigma) = \displaystyle \sum _{k = 1}^{n}\displaystyle \sum_{i,j}\mu_{k}x_{ij}\hat{u}_{ij}^{L^{\sigma_{k}}}(\sigma).$\\
$\hat{u}_{ij}^{L^{\sigma_{k}}}(\sigma)\neq 0$ for only $\sigma \in \left\{ \sigma_{1},..., \sigma_{n} \right\}$ because 
\begin{eqnarray*}
\hat{u}_{ij}^{L^{\sigma_{k}}}(\sigma)(\theta_{m},\theta_{l}) & = & \int_{G}\left < \overline{U}^{L^{\sigma}}_{t} \theta_{m},\theta_{l} \right > u_{ij}^{L^{\sigma_{k}}}(t)d\lambda(t)\\
& = & \int_{G} \overline{u}^{L^{\sigma}}_{lm}(t)u_{ij}^{L^{\sigma_{k}}}(t)d\lambda(t)\\
& = & \left\{\begin{array}{cl}
 0  &\text{  if  } \sigma \neq \sigma_{k}  \\
  & \\
\dfrac{\delta_{li}\delta_{mj}}{d_{\sigma_{k}}} &\text{ if } \sigma = \sigma_{k}.
\end{array}\right. 
\end{eqnarray*}
Hence $\hat{f} \in \mathfrak{S}_{00}(\Sigma,\mathcal{A}).$\\
Conversly let $\Phi$ be an element of $\mathfrak{S}_{00}(\Sigma,\mathcal{A})$. There exists a finite set $P = \left\{ \sigma \in \Sigma: \Phi(\sigma) \neq 0 \right\}.$ Morever each $\Phi(\sigma)$ is in the form $\Phi(\sigma) = \displaystyle\sum_{i,j=1}^{\infty} d_{\sigma}a_{ij}^{\sigma}\hat{u}_{ij}^{L^{\sigma}}.$ There exists also a finite set $I_{\sigma}$ such that $\Phi(\sigma) = \displaystyle\sum_{i,j=1}^{\infty} d_{\sigma}a_{ij}^{\sigma}\hat{u}_{ij}^{L^{\sigma}} = \displaystyle\sum_{i,j\in I_{\sigma}} d_{\sigma}a_{ij}^{\sigma}\hat{u}_{ij}^{L^{\sigma}}.$  By putting $ f = \displaystyle\sum_{\sigma \in P}d_{\sigma} \displaystyle\sum_{i,j\in I_{\sigma}} a_{ij}^{\sigma}u_{ij}^{L^{\sigma}}$ then $f\in\mathcal{L}(G) \otimes \mathcal{A}$ and we have  $\hat{f} = \displaystyle\sum_{\sigma \in P}d_{\sigma} \displaystyle\sum_{i,j\in I_{\sigma}} a_{ij}^{\sigma}\hat{u}_{ij}^{L^{\sigma}},$ thus we get  $\Phi = \hat{f}.$\\ Hence  $\widehat{\mathcal{L}(G)\otimes \mathcal{A}} =  \mathfrak{S}_{00}(\Sigma,\mathcal{A}).$ $\blacktriangledown$
\end{itemize}
\begin{lemma}\mbox{ }\\
$\mathcal{L}(G,\mathcal{A})$ is dense in $\mathfrak{C}_{0}(G,\mathcal{A})$.
\end{lemma}
\textbf{Proof}\mbox{ }\\
We identify $\mathcal{L}(G,\mathcal{A})$ with $\mathcal{L}(G,\mathbb{C})\hat{\otimes}_{\varepsilon}\mathcal{A}$ the injective tensor product of $\mathcal{L}(G,\mathbb{C})$ and $\mathcal{A}$.
According to theorem \ref{stone} \quad  $\mathcal{L}(G,\mathbb{C})$ is dense in $\mathfrak{C}_{0}(G,\mathbb{C})$ then $\mathcal{L}(G,\mathcal{A})$ is dense in $\mathfrak{C}_{0}(G,\mathcal{A})$ because $\mathfrak{C}_{0}(G,\mathcal{A})$ is isomorphic to $\mathfrak{C}_{0}(G,\mathbb{C})\hat{\otimes}_{\varepsilon}\mathcal{A}$ as a subspace of $\mathfrak{C}(G,\mathcal{A}) = \mathfrak{C}(G,\mathbb{C})\hat{\otimes}_{\varepsilon}\mathcal{A}$ according to  proposition\ref{tens1}. $\blacktriangledown$
\begin{theorem} \mbox{  }\\
 $\widehat{L_{1}(G, \mathcal{A})}$ is dense  $\mathfrak{S}_{0}(\Sigma,\mathcal{A})$ and then  for any $f$ in $L_{1}(G, \mathcal{A})$ the set $\lbrace \sigma \in \Sigma : \mbox{   } \widehat{f}(\sigma) \neq 0\rbrace$ is countable.
\end{theorem}
\textbf{Proof}\\
$\mathcal{L}(G,\mathcal{A})$ is dense in $\mathfrak{C}_{0}(G,\mathcal{A})$ for the topology of uniform convergence. Then $\mathcal{L}(G,\mathcal{A})$ is dense in $\mathcal{K}(G,\mathcal{A})$ with the same norm.
$\mathcal{K}(G,\mathcal{A})$ is also dense in $L_{1}(G,\mathcal{A})$ for $\Vert . \Vert_{1}$.\\
Let $f$ be an element of $L_{1}(G,\mathcal{A})$, $f$ is a limit of a sequence in $\mathcal{L}(G,\mathcal{A})$ then $\hat{f}$ belongs to $\overline{\widehat{\mathcal{L}(G,\mathcal{A})}} = \overline{\mathfrak{S}_{00}(\Sigma,\mathcal{A})} = \mathfrak{S}_{0}(\Sigma,\mathcal{A})$. Thus $\overline{\widehat{L_{1}(G,\mathcal{A})}} \subset \mathfrak{S}_{0}(\Sigma,A)$. Conversly $\mathfrak{S}_{0}(\Sigma,\mathcal{A}) = \overline{\mathfrak{S}_{00}(\Sigma,\mathcal{A})} = \overline{\widehat{\mathcal{L}(G,\mathcal{A})}} \subset \overline{\widehat{L_{1}(G,\mathcal{A})}}$ It implies $ \mathfrak{S}_{0}(\Sigma,A) \subset \overline{\widehat{L_{1}(G,\mathcal{A})}}$. Hence  $\overline{\widehat{L_{1}(G,\mathcal{A})}} = \mathfrak{S}_{0}(\Sigma,\mathcal{A}).$\\
Morever $\lbrace \sigma \in \Sigma : \mbox{   } \widehat{f}(\sigma) \neq 0\rbrace$ is countable because for every $n \in \mathbb{N}^{*}$ the set $\lbrace  \sigma \in \Sigma$: 
 $\Vert \widehat{f}(\sigma)  \Vert > \dfrac{1}{n} \rbrace$ is finite. $\blacktriangledown$
 \subsection{Plancherel theorem}
 \begin{proposition} \label{Scalar}
 Suppose that $\mathcal{A}$ is a Hilbert space then the mapping
 \begin{eqnarray}
(\Phi,\Psi) & \longmapsto & \left < \Phi,\Psi \right > := \sum_{\sigma\in\Sigma}d_{\sigma}\sum_{i,j = 1}^{\infty}\left < \Phi(\sigma)(\theta_{i},\theta_{j}) , \Psi(\sigma)(\theta_{i},\theta_{j}) \right >_{\mathcal{A}}
 \end{eqnarray}
 is a scalar product in $\mathfrak{S}_{2}(\Sigma, \mathcal{A})$.
 \end{proposition}
 \textbf{Proof}
 \begin{eqnarray*}
 \left\vert \sum_{\sigma\in\Sigma}d_{\sigma}\sum_{i,j = 1}^{\infty}\left < \Phi(\sigma)(\theta_{i},\theta_{j}) , \Psi(\sigma)(\theta_{i},\theta_{j}) \right >_{\mathcal{A}} \right\vert  & \leq & \sum_{\sigma\in\Sigma}\sum_{i,j = 1}^{\infty}d_{\sigma}\left \vert\left < \Phi(\sigma)(\theta_{i},\theta_{j}) , \Psi(\sigma)(\theta_{i},\theta_{j}) \right >_{\mathcal{A}}\right \vert\\
& \leq & \sum_{\sigma\in\Sigma}\sum_{i,j = 1}^{\infty}d_{\sigma}^{\frac{1}{2}}\left \Vert \Phi(\sigma)(\theta_{i},\theta_{j}) \right \Vert_{\mathcal{A}}d_{\sigma}^{\frac{1}{2}}\left \Vert \Psi(\sigma)(\theta_{i},\theta_{j}) \right \Vert_{\mathcal{A}}\\
& \leq &  \sum_{\sigma\in\Sigma}\sum_{i,j = 1}^{\infty}\left(d_{\sigma}\left \Vert \Phi(\sigma)(\theta_{i},\theta_{j}) \right \Vert_{\mathcal{A}}^{2}\right)^{\frac{1}{2}} \times\\
&  & \sum_{\sigma\in\Sigma}\sum_{i,j = 1}^{\infty}\left(d_{\sigma}\left \Vert \Psi(\sigma)(\theta_{i},\theta_{j}) \right \Vert_{\mathcal{A}}^{2}\right)^{\frac{1}{2}}\\
& < & \infty.
 \end{eqnarray*}
 Then the mapping is well defined. The sequel of the proof is trivial, it depends to the fact that $\left < , \right >_{\mathcal{A}}$ is a scalar product in  $\mathcal{A}.$  $\blacktriangledown$\\

The following theorem is the Plancherel theorem in this context. 
 
\begin{theorem}
Let $\mathcal{A}$ be a Hilbert space then the mapping $f \longmapsto \hat{f}$ is an isometry of $L_{2}(G, \mathcal{A}, \lambda)$ into  $\mathfrak{S}_{2}(\Sigma, \mathcal{A})$ and then $\mathfrak{S}_{2}(\Sigma, \mathcal{A})$ is a Hilbert space.
\end{theorem}
\textbf{Proof}
Since $\mathcal{A}$ is a Hilbert space then so is $L_{2}(G, \mathcal{A}, \lambda)$.

Let $f$ be an element of $L_{2}(G, \mathcal{A}, \lambda)$, put $a_{ij}^{\sigma} = \hat{f}(\sigma)(\theta_{i},\theta_{j})$
\begin{eqnarray*}
\Vert f \Vert_{2}^{2} & = & \left <\sum_{\sigma\in \Sigma}\sum_{i,j}^{\infty}d_{\sigma}a_{ij}^{\sigma}u_{ij}^{L^{\sigma}}, \sum_{\sigma\in \Sigma}\sum_{i,j}^{\infty}d_{\sigma}a_{ij}^{\sigma}u_{ij}^{L^{\sigma}} \right > \\
& = & \sum_{\sigma\in \Sigma}\sum_{i,j}^{\infty}d_{\sigma}^{2}\Vert a_{ij}^{\sigma}\Vert_{\mathcal{A}}^{2}\Vert u_{ij}^{L^{\sigma}} \Vert_{2}^{2}\\
& = & \sum_{\sigma\in \Sigma}\sum_{i,j}^{\infty}d_{\sigma}\Vert a_{ij}^{\sigma}\Vert_{\mathcal{A}}^{2}
\end{eqnarray*}
because $d_{\sigma}\Vert u_{ij}^{L^{\sigma}} \Vert_{2}^{2} = 1$ according to Shur orthogonality theorem.

Now $\Vert \hat{f} \Vert ^{2} = \displaystyle\sum_{\sigma\in \Sigma}\sum_{i,j}^{\infty}d_{\sigma}\Vert a_{ij}^{\sigma}\Vert_{\mathcal{A}}^{2}$. We have $\Vert \hat{f} \Vert ^{2} = \displaystyle\sum_{\sigma\in \Sigma}\sum_{i,j}^{\infty}d_{\sigma}\Vert a_{ij}^{\sigma}\Vert_{\mathcal{A}}^{2}  = \Vert f \Vert ^{2} < \infty$ then $\hat{f}\in \mathfrak{S}_{2}(\Sigma, \mathcal{A})$.

Conversly 
 \section*{Concluding remarks}
The aspects approached in this paper are consequences of the work on the Fourier-Stieltjes transform on locally compact group recently elaborated in \cite{yao}.

We have through this work given meaning to the theory developed in \cite{yao}. We have shown that the Fourier-Stieltjes transform defined in a locally compact group satisfies the main and classic theorems in analysis. We have among other established the theorems of Lebesgue and that of Plancherel then the identity of Parseval relating to our concept.

We intend in the future to study, by way of application, this theory in certain
much more particular groups such as the linear groups $GL(n,\mathbb{C})$ and some of its subgroups.

\end{document}